\documentclass[11pt,a4paper]{article}
\usepackage{fontenc,amsmath,mathrsfs,amsfonts,amsthm, verbatim,geroENG,amssymb, MnSymbol, enumerate, multicol, color, hyperref}
\usepackage{pgf,tikz}
\usepackage{mathrsfs}
\usepackage[english]{babel}
\usepackage[top=1.25in, bottom=1.25in, left=1.0in, right=1.0in]{geometry}
\title{Good's Theorem for Hurwitz Continued Fractions}
\author{Gerardo Gonz\'alez Robert}
\begin{document}
\date{}
\maketitle

\begin{abstract}
Good's Theorem for regular continued fraction states that the set of real numbers $[a_0;a_1,a_2,\ldots]$ such that $\displaystyle\lim_{n\to\infty} a_n=\infty$ has Hausdorff dimension $\tfrac{1}{2}$. We show an analogous result for the complex plane and Hurwitz Continued Fractions: the set of complex numbers whose Hurwitz Continued fraction $[a_0;a_1,a_2,\ldots]$ satisfies $\displaystyle\lim_{n\to\infty} |a_n|=\infty$ has Hausdorff dimension $1$, half of the ambient space's dimension.
\end{abstract}

\section{Introduction}
A classical result on Diophantine approximation states that a real number $\alpha$ is irrational if and only if there are infinitely many co-prime $p,q\in\Za$, $q\geq 1$, for which
\[
\left| \alpha - \frac{p}{q}\right| < \frac{1}{q^2}.
\]
For some numbers $\alpha$, the exponent of $q$ in the previous inequality cannot be improved. Such $\alpha$ are called badly approximable and the set of badly approximable real numbers is denoted by $\bad_{\RE}$. More precisely, $\alpha\in\RE\setminus \QU$ is \textbf{badly approximable} if there exists $c>0$ such that
\[
\forall (p,q)\in\Za\times \Na\quad
\left| \alpha - \frac{p}{q}\right| \geq  \frac{c}{q^2}.
\]
It is well known that $\bad_{\RE}$ is exactly the set of irrational numbers whose regular continued fraction is given by a bounded sequence. This characterization allows us to show by elementary means that $\leb_1\left(\bad_{\RE}\right)=0$ (see \cite{khin}, Theorem 23), where $\leb_1$ denotes the Lebesgue measure in $\RE$. Although small in terms of Lebesgue measure, V. Jarn\'ik showed in 1928 (Satz 4, \cite{jarnik28}) that $\bad_{\RE}$ is rather large.

\begin{teo01}[V. Jarn\'ik,1928]
The set $\bad_{\RE}$ has full Hausdorff dimension; that is
\[
\dim_H \left\{[0;a_1,a_2,\ldots] \in \RE: \limsup_{n\to\infty} a_n<+\infty\right\} = 1.
\]
\end{teo01}

On the basis of Jarn\'ik's techniques, I.J. Good computed in 1941 the Hausdorff dimension of similar sets (Theorem 1, \cite{good}). 
\begin{teo01}[I. G. Good, 1941]\label{TeoGood1941}
The following equality holds
\[
\dim_H \left\{[0;a_1,a_2,\ldots] \in \RE: \lim_{n\to\infty} a_n = +\infty\right\} = \frac{1}{2}.
\]
\end{teo01}

Several extensions of the continued fraction theory and of Jarn\'ik's Theorem  have been successfully carried out (v. gr. \cite{ktv} or \cite{schweiger} and the references therein). Recently, S.G. Dani and A. Nogueira developed in \cite{daninog} complex continued fractions algorithms which include the one suggested by Adolf Hurwitz in \cite{hur87}. We describe Hurwitz continued fractions in detail in Section 2; in the mean time, it suffices to know that Hurwitz continued fractions associate to each irrational $\zeta$ a sequence of non-zero Gaussian integers $\sanu$ such that
\[
\zeta=[a_0;a_1,a_2,\ldots]:= a_0 + \cfrac{1}{a_1 + \cfrac{1}{a_2+ \cfrac{1}{\ddots}}}
\]
(understood in the usual sense).

As in the real case, by Dirichlet's Pigeonhole Principle, there is an absolute constant $C>0$ such that a complex number $\zeta$ is irrational, i.e. $\zeta\in\Cx\setminus\QU(i)$, if and only if there are infinitely many co-prime $p,q\in\Za[i]$, $|q|\geq 1$, satisfying
\[
\left| \zeta - \frac{p}{q}\right| \leq \frac{C}{|q|^2}.
\]
We say that a complex irrational $\zeta$ is \textbf{badly approximable} if, for $\Za[i]^*:=\Za[i]\setminus\{0\}$, 
\[
\exists c>0 \quad
\forall (p,q)\in\Za[i]\times\Za[i]^* \quad
\left| \zeta - \frac{p}{q}\right| \geq \frac{c}{|q|^2},
\]
and we denote by $\bad_{\Cx}$ the set of badly approximable complex numbers. 

Some properties of $\bad_{\Cx}$ are well known. For example, $\bad_{\Cx}$ is null with respect to the Lebesgue measure, it is $\tfrac{1}{2}$-winning in the sense of Schmidt games and, hence, it has full Hausdorff dimension (see Theorem 5.2, \cite{dodkrist}). $\bad_{\Cx}$ can also be characterized in terms of Hurwitz Continued Fractions: a complex irrational $\zeta$ belongs to $\bad_{\Cx}$ if and only if its Hurwitz Continued fraction is bounded (see Theorem 1 in \cite{hines} or Theorem 4.1 in \cite{gero18}). Our main result, Theorem \ref{TeoGero}, is also an analogy between regular and Hurwitz continued fractions.

\begin{teo01}\label{TeoGero}
For any $z\in\Cx\setminus\QU(i)$, let $z=[a_0;a_1,a_2,\ldots]$ denote its Hurwitz continued fraction. The following equality holds
\[
\dim_H \left\{ z=[a_0;a_1,a_2,a_3,\ldots] \in\Cx :  \lim_{n\to\infty} |a_n|=+\infty\right\}  =1.
\]
\end{teo01}
Although similarities between real and complex continued fractions abound, some differences have to be considered. For example, while regular continued fractions establish an homeomorphism between the irrationals in $[0,1]$ and $\Na^{\Na}$, the space of sequences associated to Hurwitz continued fractions is much more complicated. Since the difficulty arises from sequences $\sanu$ with $\min_{n\geq 1} |a_n|\leq \sqrt{8}$, it is natural to ask how large are the subsets of $\bad_{\Cx}$ where the absolute value of the terms of the Hurwitz continued fraction satisfy a uniform lower bound. In this direction, the proof of Theorem \ref{TeoGero} gives us Corollary \ref{CoroGero} (see Section 2 for the definition of $\mfF$).

For any $L>0$ define
\[
E_L:=\left\{z=[0;a_1,a_2,a_3,\ldots] \in\mfF : \forall n\in\Na \quad L\leq |a_n|\right\}.
\]

\begin{coro01}\label{CoroGero}
$
\displaystyle \lim_{L\to\infty} \dim_H  E_{L}\cap \bad_{\Cx} = 1.
$
\end{coro01}

The organization of the text is as follows. In Section 2, we define precisely the Hurwitz continued fraction algorithm and discuss some of its properties. In Section 3, we give two lemmas for estimating the Hausdorff dimension of a class of Cantor sets in complete metric spaces. In Section 4 we state some preliminary lemmas concerning the Hausdorff dimension of sets obtained by imposing restrictions on the Hurwitz continued fraction expansions. In Section 5 we show Theorem \ref{TeoGero} and Corollary \ref{CoroGero}. Finally, in Section 6 we prove the preliminary lemmas stated in Section 4.
\paragraph{Notation.}
\begin{enumerate}[(1)]
\item $\Na$ is the set of natural numbers, considered as the set of positive integers.
\item If $(X,d)$ is a metric space and $A\subseteq X$, the diameter of $A$ is  $|A|:=\sup\{d(x,y):x,y\in A\}$. 
\item For any complex number $z$ we write $\Dx(z)=\{w\in\Cx: |z-w|< 1\}$, $\overline{\Dx}(z)$ is the closure of $\Dx(z)$ and $C(z)$ is its boundary.
\item Let $A \subseteq \Cx$. $A^{\circ}$ is the interior of $A$, $\Cl(A)$ is the closure of $A$, and $A^{-1}=\{z^{-1}: z\in A\}$.
\item For $A\subseteq \Cx$, $\#A$ is the number of elements contained in $A$. When $A$ is infinite, we just write $\#A=\infty$.
\end{enumerate}

\section{Hurwitz Continued Fractions}
Denote by $[\cdot]_{\RE}:\RE\to\Za$ the function that assigns to each $x\in\RE$ the largest integer less than or equal to $x$. Let $[\cdot]:\Cx\to\Za[i]$ be given by
\[
 \forall z\in\Cx \quad [z]=\left[ \Re z + \frac{1}{2} \right]_{\RE} + i\left[ \Im z+ \frac{1}{2} \right]_{\RE}.
\]
Denote by $\mfF$ the inverse image of $0$ under $[\cdot]$,
\begin{equation}\label{EcDefmfF}
\mfF:=\left\{ z\in\Cx: -\frac{1}{2}\leq \Re z, \Im z<\frac{1}{2}\right\}, \quad \text{and define}\quad \mfF^*:=\mfF\setminus\{0\}.
\end{equation}
Let $T:\mfF^*\to\mfF$, with $\mfF^*:=\mfF\setminus\{0\}$, be given by
\[
\forall z\in\mfF^* \quad T(z) = \frac{1}{z} - \left[ \frac{1}{z}\right].
\]
For any $z\in\mfF^*$ define---as long as the operations make sense---the sequences $(a_n)_{n\geq 0}$, $(z_n)_{n\geq 1}$ by
\begin{align*}
z_1:=z, &\quad
\forall n\in\Na \quad z_{n+1}=T(z_{n}), \nonumber\\
a_0:= 0, 
&\quad 
\forall n\in\Na \quad 
a_{n} = \left[z_n^{-1}\right]. \nonumber
\end{align*}
The \textbf{Hurwitz continued fraction} (HCF) of $z$ is the sequence $\sanu$. We can easily extend the definition to an arbitrary complex number. Indeed, the Hurwitz continued fraction of $w\in \Cx$ is the sequence $\sean$ where $a_0=[w]$ and $\sanu$ is the HCF of $w-a_0$. We refer to the numbers $a_n$ as \textbf{elements}. Following \cite{daninog}, we call the sequences $\sepn$, $\seqn$ given by
\[
\begin{pmatrix}
p_{-1} & p_0 \\
q_{-1} & q_0 \\
\end{pmatrix}
=
\begin{pmatrix}
1 & 0 \\
0 & 1 \\
\end{pmatrix},
\quad\quad
\forall n\in\Na \quad
\begin{pmatrix}
p_{n}  \\
q_{n} \\
\end{pmatrix}
=\begin{pmatrix}
p_{n-1} & p_{n-2} \\
q_{n-1} & q_{n-2} \\
\end{pmatrix}
\begin{pmatrix}
a_{n} \\
1\\
\end{pmatrix}.
\]
the \textbf{$\mathcal{Q}$-pair} of $z\in\mfF$. 

We summarize some well-known properties of HCF.
\begin{propo01}\label{Propo2.1}
Let $z\neq 0$ belong to $\mfF$ and let $\sean$, $\sepn$, $\seqn$ be its associated sequences.
\begin{enumerate}[i.]
\item \label{Propo2.1_i}For every $n\in\Na$
\[
\frac{p_n}{q_n} = [0;a_1,a_2,\ldots,a_n]:=\cfrac{1}{a_1+ \cfrac{1}{a_2+\cfrac{1}{\ddots + \cfrac{1}{a_n}}}}.
\]
\item \label{Propo2.1_ii} The sequence $\sean$ is infinite if and only if $z\in\mfF\setminus\QU(i)$. In this case, we have that
\[
\lim_{n\to\infty} [0;a_1,a_2,\ldots,a_n] = z.
\]
Thus, HCF give an injection from $\mfF\setminus\QU(i)$ into $\Za[i]^{\Na}$.
\item \label{Propo2.1_iii} The sequence $(|q_n|)_{n\geq 0}$ is strictly increasing. Moreover, there exists a number $\psi>1$ such that $|q_n|\geq \psi^{n-1}$ for all $n\in\Na_0$.
\item \label{Propo2.1_iv} If $z=[0;a_1,a_2,\ldots]\in\mfF\setminus\QU(i)$, we have
\[
\forall n\in\Na \quad z= \frac{(a_{n+1}+[0;a_{n+2},a_{n+3},\ldots])p_{n}+p_{n-1}}{(a_{n+1}+[0;a_{n+2},a_{n+3},\ldots])q_{n}+q_{n-1}}.
\]
\item \label{Propo2.1_v} For every $n\in\Na$ such that $p_n$ and $q_n$ are defined
\[
\left| z - \frac{p_n}{q_n}\right| \leq \frac{1}{|q_n|^2}.
\]
\end{enumerate}
\end{propo01}
\begin{proof}
Part \ref{Propo2.1_i}. is trivial, \ref{Propo2.1_ii}. is Theorem 6.1. in \cite{daninog}, the monotonicity in \ref{Propo2.1_iii}. is on page 195 of \cite{hur87} and the exponential growth is Corollary 5.3. in \cite{daninog}, \ref{Propo2.1_iv}. is a restatement of Proposition 3.3. in \cite{daninog}, and \ref{Propo2.1_v}. is Theorem 1 in \cite{lakein01}. 
\end{proof}

Define $I=\{a\in\Za[i]: |a|\geq \sqrt{2}\}$. For every $a\in \Za[i]$ the \textbf{cylinder of level $1$} is
\[
\clC_1(a) = \left\{z\in \mfF: [z^{-1}] = a \right\}.
\]
Note that $\clC_1(a)\neq \vac$ if and only if $a\in I$ and that $\{\clC_1(a)\}_{a\in I}$ is a countable partition of $\mfF$. We can define the HCF elements as functions from $\mfF$ onto $I$ as follows
\[
\forall n\in \Na \quad\qquad a_n(z) = a \;\iff\; T^{n-1}(z) \in \clC_1(a),
\]
where $T^0:\mfF\to\mfF$ is the identity and $T^{n}=T\circ T^{n-1}$ for $n\in\Na$. For $\bfa=(a_1,\ldots,a_n)\in I^n$ or $\bfa\in I^{\Na}$ the \textbf{cylinder of level} $n\in\Na$, $\clC_n(\bfa)$, is the set
\[
\clC_n(\bfa)=\left\{ z\in\mfF: a_1(z)=a_1,\ldots, a_n(z)= a_n\right\}.
\]
A sequence $\bfa\in I^{\Na}$ is \textbf{admissible} or \textbf{valid} if $\clC_n(\bfa)\neq \vac$ for every $n\in\Na$. We denote the set of admissible sequences by $\Omega^{\HCF}_0$ and define $\Phi:\Omega^{\HCF}_0\to \mfF$ by
\[
\forall \bfa\in\Omega_0^{\HCF} 
\quad \Phi(\bfa)=[0;a_1,a_2,a_3,\ldots].
\]
A \textbf{maximal feasible set} is a set $\widetilde{\mfF}\subseteq \mfF$ such that for some $\bfa\in\Omega_0^{\HCF}$ and some $n\in\Na$ we have
\[
T^n\left[ \clC_n(\bfa) \right] =\widetilde{\mfF}.
\]
Maximal feasible sets may be proper subsets of $\mfF$. For instance, $w\in T_1[\clC_1(-2)]$ if and only if $(-2+w)^{-1}\in\mfF$, which is equivalent to
\[
w\in (2+\mfF^{-1}) \cap \mfF = \mfF\setminus \Dx(1).
\]
Hence, $T[\clC_1(-2)]=\mfF\setminus \Dx(1)$. In a like manner, we can show that $T^2[\clC_2(-2,1-3i)]=\{it: -\tfrac{1}{2}\leq t< \tfrac{1}{2}\}$. We can also prove similarly that 
$T[\clC_1(a)]= \mfF$ if and only if $|a|\geq \sqrt{8}$. Moreover, when we dismiss the boundary, all the non-empty sets $T[\clC_1(a)]$ are of the form $i^j\mfF_k$ for some $j,k\in \{1,2,3,4\}$, where
\begin{align*}
\mfF_1 := \left(\mfF\setminus \left( \Dx(-1) \cup  \Dx(-i)\right)\right)^{\circ}, 
&\quad
\mfF_2 := \left(\mfF\setminus \Dx(-1)\right)^{\circ},
\nonumber\\
\mfF_3 := \left(\mfF\setminus \Dx(-1-i)\right)^{\circ},
&\quad
\mfF_4=\mfF^{\circ}. \nonumber
\end{align*}
By computing directly the inversion of some circles (v.gr. $C(1+i)^{-1}=C(1-i)$), it can be shown inductively that these thirteen sets exhaust all the possibilities for the shapes that the maximal feasible sets with non-empty interior may assume (cfr. \cite{hiva18}).

\begin{rema1000}
The proof of Theorem \ref{TeoGero} only uses maximal feasible sets with non-empty interior.
\end{rema1000}

\subsection{Some sets defined via HCF}
Let us start with an observation.
\begin{lem01}\label{Le2.0}
There exists an absolute constant $k>0$ such that any for $\alpha,\beta\in\mfF$ with $|a_1(\alpha)|\geq \sqrt{8}$, $|a_1(\beta)|\geq \sqrt{8}$ and any $a,b\in\Za[i]$, $a\neq b$,
\[
|(a+\alpha) - (b+\beta)|>k.
\]
\end{lem01}
\begin{proof}
The result follows from the triangle inequality.
\end{proof}
Recall that for any $L>0$   
\[
E_{L}:=\left\{ [0;a_1,a_2,\ldots]\in\mfF: \forall n\in \Na \quad |a_n|\geq L \right\}.
\]
An inductive argument tells us that if $L\geq 8$, then every $z=[0;a_1,a_2,\ldots]\in E_L$ satisfies $T^n[\clC_n(a_1,\ldots,a_n)]=\mfF$ for all $n$. Hence, every sequence $\sanu\in\Za[i]^{\Na}$ verifying $|a_n|\geq \sqrt{8}$, $n\in\Na$, belongs to $\Omega_0^{\HCF}$.

\begin{lem01}\label{Le2.1}
There is a constant $\gamma>0$ such that 
\begin{equation}\label{Ec2.2}
\forall \bfa\in E_{\sqrt{8}} \quad \forall n\in\Na \quad 
\frac{\gamma}{|q_n(\bfa)|^2} \leq |\clC_n(\bfa)\cap E_{\sqrt{8}}| \leq |\clC_n(\bfa)|\leq \frac{2}{|q_n(\bfa)|^2}.
\end{equation}
\end{lem01} 
\begin{proof}
Let $\bfa,n$ be as in the statement. The inequality $|\clC_n(\bfa)\cap E_{\sqrt{8}}|\leq |\clC_n(\bfa)|$ is trivial and $|\clC_n(\bfa)|\leq 2|q_n|^{-2}$ follows from Proposition \ref{Propo2.1}. In order to show the remaining inequality, we estimate the distance between two particular elements of $\clC_n(\bfa)\cap E_{\sqrt{8}}$. Define $\xi=[3+4i;3+4i,3+4i,\ldots]$ (a constant sequence). Let $z,w\in\clC_n(\bfa)\cap E_{\sqrt{8}}$ be given by
\[
z=[0;a_1,\ldots,a_n,3+4i,3+4i,3+4i,\ldots], \quad
w=[0;a_1,\ldots,a_n,-3-4i,-3-4i,-3-4i,\ldots].
\]
Then, if $\sepn$, $\seqn$ are the $\clQ$-pair of $\bfa$, by Proposition \ref{Propo2.1} we have
\begin{align*}
|\clC_n(\bfa)\cap E_{\sqrt{8}}| &\geq |z-w| \nonumber\\
&\geq \left| \frac{p_n\xi + p_{n-1}}{q_n\xi + q_{n-1}} - \frac{p_n(-\xi) + p_{n-1}}{q_n(-\xi) + q_{n-1}} \right| = \frac{1}{|q_n|^2} \frac{2|\xi|}{\left| \left( \xi + \frac{q_{n-1}}{q_n}\right)\left( -\xi + \frac{q_{n-1}}{q_n}\right)\right|} \nonumber\\
& \geq
\frac{1}{|q_n|^2}\frac{2|\xi|}{(|\xi|+1)^2}. \nonumber
\end{align*}
\end{proof}

For any $\bfa\in\Omega_0^{\HCF}$, define the transformations $(t_{\bfa,n})_{n\geq 1}$ and $(v_{\bfa,n})_{n\geq 1}$ by
\[
\forall n\in\Na \quad \forall z\in\Cx\quad
t_{\bfa,n}(z) = \frac{q_nz-p_n}{-q_{n-1}z+p_{n-1}},\quad
v_{\bfa,n}(z) = \frac{p_{n-1}z+p_n}{q_{n-1}z+q_n},
\]
where $\sepn$, $\seqn$ are the $\clQ$-pair of $\bfa$. Note that for every $n\in\Na$ the functions $t_{\bfa,n}$ and $v_{\bfa,n}$ are inverses, and that the restrictions of $T^n$ and $t_{\bfa,n}$ to $\clC_n(\bfa)$ coincide. Therefore, for every $n$ the restriction of $T^n$ to $\clC_n(\bfa)$ is bi-Lipschitz onto its image.

\begin{lem01}\label{Le2.2}
For any $\bfa\in\Omega^{\HCF}_0$ and any $n\in\Na$, the map $t_{\bfa,n}:\clC_n(\bfa)\to T^n[ \clC_n(\bfa)]$ acts via
\[
[0;a_1,\ldots,a_n,x_1,x_2,\ldots] \mapsto [0;x_1,x_2,\ldots],
\]
and is bi-Lipschitz. Moreover, if $\bfb\in\Omega^{\HCF}_0$, $m\in\Na$, and $T^n[\clC_n(\bfa)]=T^m[\clC_m(\bfb)]$, then the map from $\clC_n(\bfa)$ to $\clC_m(\bfb)$ given by
\[
[0;a_1,\ldots,a_n,x_1,x_2,\ldots] \mapsto [0;b_1,\ldots,b_m,x_1,x_2,\ldots]
\]
is bi-Lipschitz. 
\end{lem01}
Let us define for any $L>0$ and any $N\in\Na$  
\begin{align}\label{DefE_I(N)}
E_L(N) &:=\left\{z=[0;a_1,a_2,\ldots] \in\mfF: \; \forall n\in\Na_0 \quad |a_{N+n}|\geq L\right\}, \nonumber\\
E_L' &:=\left\{z=[0;a_1,a_2,\ldots] \in\mfF: \; \liminf_{n\to\infty} |a_{n}|\geq L\right\}. \nonumber
\end{align}
so  $E_L(1)=E_L$ for all $L>0$.
\begin{lem01}\label{LeDesc01}
For any $L>0$, we have $\dim_H E_L = \dim_H E_L'$.
\end{lem01}
\begin{proof}
Let $L>0$. By definition, $E_L'=\bigcup_{N\geq 1} E_L(N)$, hence
\[
\dim_H E_L' = \sup_{n\in\Na} \dim_H E_L(N).
\]
The result will follow if we show that $(\dim_H E_L(N))_{N\geq 1}$ is constant. On the one hand, $E_L\subseteq E_L(N)$ for every $N\in\Na$ and hence 
\[
\forall n\in\Na\quad \dim_H E_L \leq \dim_H E_L(N). 
\]
On the other hand, take $N\in\Na_{\geq 2}$. Define $G_L:=\{a\in\Za[i]:|a|\geq L\}$, then $E_L(N)$ can be written as a countable union
\[
E_L(N)= \bigcup_{\bfa\in I^{N-1}} \Phi\left[ \left(\{\bfa\}\times G_L^{\Na}\right)\cap \Omega_0^{\HCF}\right]
\]
(some terms might be empty), therefore
\[
\dim_H E_L(N) = \sup_{\bfa\in I^{N-1}} \dim_H \Phi\left[ \left(\{\bfa\}\times G_L^{\Na}\right)\cap \Omega_0^{\HCF}\right].
\]
Take $\bfa\in I^{N-1}$ such that $\clC_{N-1}(\bfa)\neq \vac$. Lemma \ref{Le2.2} tells us that the map
\begin{align*}
\Phi\left[ \left(\{\bfa\}\times G_L^{\Na}\right)\cap \Omega_0^{\HCF}\right] &\to T^{N-1}\left[ \clC_{N-1}(\bfa)\right] \cap E_L, \nonumber\\
[0;a_1,\ldots,a_n,x_1,x_2,\ldots] &\mapsto [0;x_1,x_2,\ldots]\nonumber
\end{align*}
is bi-Lipschitz. The invariance of the Hausdorff dimension under bi-Lipschitz maps gives
\[
\dim_H \Phi\left[\left(\{\bfa\}\times G_L^{\Na} \cap \Omega_0^{\HCF}\right]\right)
=
\dim_H T^{N-1}\left[ \clC_{N-1}(\bfa)\right] \cap E_L 
\leq 
\dim_H E_L.
\]
Taking the supremum over $\bfa\in I^{N-1}$, we conclude $\dim_H E_L(N)\leq \dim_H E_L$.
\end{proof}

\section{Generalized Jarn\'ik Lemmas}
Two results, which we shall call Generalized Jarn\'ik Lemmas, lie in the heart of Theorem \ref{TeoGero}. In this section, our framework is a more restricted version of strongly tree-like sets as defined in \cite{klewei10}.
\begin{def01}\label{DefDSTL}
Let $(X,d)$ be a complete metric space. A family of compact sets $\clA$ is \textbf{diametrically strongly tree-like} if $\clA=\bigcup_{n=0}^{\infty} \clA_n$ where each $\clA_k$ is finite, $\#\clA_0 = 1$, and
\begin{enumerate}[i.]
\item $\forall A\in\clA \quad |A|>0$,
\item $\forall n\in\Na \quad \forall A,B\in\clA_n \quad (A=B)\lor (A\cap B=\vac)$,
\item $\forall n\in\Na \quad \forall B\in\clA_n \quad \exists A\in\clA_{n-1}\quad B\subseteq A$,
\item $\forall n\in\Na \quad \forall A\in\clA_{n-1} \quad \exists B\in\clA_{n}\quad B\subseteq A$,
\item $d_n(\clA):=\max\{|A|:A\in\clA_n\}\to 0$ as $n\to\infty$.
\end{enumerate}
For each $n\in\Na_0$ and each $A\in\clA_n$, the \textbf{descendants} of $A$ are the members of $D(A)=\{B\in\clA_{n+1}: B\subseteq A\}$. The quantity $d_n(\clA)$ is the $n$\textbf{-th stage diameter}. The \textbf{limit set} of $\clA$, $\mathbf{A}_{\infty}$, is 
\[
\mathbf{A}_{\infty}:= \bigcap_{n=0}^{\infty} \bigcup_{A\in\clA_n} A.
\]
\end{def01}

As an example, for $j\in\{0,2\}$ consider $f_j:[0,1]\to[0,1]$, $f(x)=\tfrac{x}{3}+\tfrac{j}{2}$. Then, the family of compact sets $\clA=\bigcup_{n\geq 0} \clA_n$ given by
\[
\clA_0:=\{[0,1]\}, \qquad \forall n\in\Na \quad \clA_n:=\{f_{a_1}\circ\ldots\circ f_{a_n}[0,1]:a_1,\ldots,a_n\in\{0,2\}\}
\]
is diametrically strongly tree-like and the corresponding $\bfA_{\infty}$ is the middle third Cantor set.

As in \cite{jarnik28}, if $\mfY=\{Y_j\}_{j\in J}$ is an at most countable family of subsets of $X$ and $s\geq 0$, we write
\[
\Lambda_s(\mfY):= \sum_{j\in J} |Y_j|^s.
\]
For $X_1,X_2\subseteq X$ we denote the distance between them by
\[
d(X_1,X_2):=\inf\{ d(x_1,x_2): x_1\in X_1, x_2\in X_2\}.
\]

\begin{lem01}[First Generalized Jarn\'ik Lemma]\label{LeGJL01}
Let $\clA$ be a diametrically strongly tree-like family of compact sets with limit set $\bfA_{\infty}$. Suppose that
\begin{equation}\label{Ec1GJL01}
\liminf_{n\to\infty} \frac{\log\left( d_n(\clA)^{-1}\right)}{n}>0,
\end{equation}
and that there exists a sequence in $(0,1)$, $(B_n)_{n\geq 1}$, such that
\begin{equation}\label{Ec1GJL02}
\limsup_{n\to\infty} \frac{\log\log\left(B_n^{-1}\right)}{\log n}<1,
\end{equation}
\begin{equation}\label{Ec1GJL03}
\forall n\in\Na \quad \forall A\in\clA_n \quad \forall\,Y,Z\in D(A) 
\quad 
\left(Y\neq Z \implies d(Y,Z)\geq B_n |A|\right).
\end{equation}
If for $s>0$ there exists some $c>0$ satisfying
\begin{equation}\label{Ec1GJL04}
\forall \mfX\in 2^{\clA} \quad \left( \#\mfX<+\infty  \quad\&\quad  \mathbf{A}_{\infty}\subseteq \bigcup_{A\in\mfX} A  
\quad\implies\quad
\Lambda_s\mfX>c\right),
\end{equation}
then $\dim_H A_{\infty}\geq s$.
\end{lem01}

\begin{proof}
Let us keep the statement's notation. We divide the argument into two parts.

First, let $\widetilde{G}$ be an open set such that $\widetilde{G}\cap\bfA_{\infty}\neq\vac$ and $|\widetilde{G}|<+\infty$. Define
\[
\widetilde{G}':= \widetilde{G}\cap \bfA_{\infty}, \quad
n(\widetilde{G}):=\max\{m\in\Na_0: \exists A\in \clA_m \quad \widetilde{G}'\subseteq A\},
\]
and let $A(\widetilde{G})\in \clA_{n(\widetilde{G})}$ satisfy $\widetilde{G}'\subseteq A(\widetilde{G})$. By definition of $n(\widetilde{G})$, there are $Y,Z\in D(A({\widetilde{G}}))$ such that $Y\neq Z$, $Y\cap \widetilde{G}\neq\vac$, and $Z\cap \widetilde{G}\neq\vac$; hence,
\begin{equation}\label{EcPr1GJL01}
\min\{|\widetilde{G}|, |A(\widetilde{G})|\}\geq |\widetilde{G}'| \geq B_{n(\widetilde{G})} |A(\widetilde{G})|.
\end{equation}
Since each $\clA_k$ is finite,
\[
\forall m\in\Na\quad \alpha_m:=\min\left\{|A|: A\in \bigcup_{k=1}^m \clA_k\right\}>0,
\]
and $(B_m\alpha_m)_{m\geq 1}$ consists of strictly positive terms. Therefore, by \eqref{EcPr1GJL01},  $n(\widetilde{G})\to\infty$ when $|\widetilde{G}|\to 0$. Note also that \eqref{Ec1GJL01} implies the existence of some $\kappa>1$ such that 
\begin{equation}\label{EcPr1GJL02}
\text{ for every large } j\in\Na \qquad d_j(\clA)\leq \frac{1}{\kappa^j}.
\end{equation}
For the second part, let $0<\veps<s$. Take $\delta>0$ and let $\clG$ be a finite open cover of $\bfA_{\infty}$ with $|G|<\delta$ for all $G\in\clG$. By \eqref{Ec1GJL02} and the first part of the proof, we can pick $\delta>0$ so small that
\[
\forall j\in\Na 
\qquad 
j\geq \min\{ n(G): G\in\clG\} \implies \kappa^{\veps j}B_j^s>1.
\]
Let $A\in\clA$ be such that $A_{\clG}:=\{G\in\clG:A(G)=A\}\neq\vac$ and let $j_0\in\Na$ satisfy $A\in\clA_{j_0}$. It follows from \eqref{EcPr1GJL01} and \eqref{EcPr1GJL02} that 
\begin{align*}
|A|^s  &\leq \frac{1}{B_{j_0}^s}\sum_{G\in A_{\clG}} |G'|^s = \frac{1}{B_{j_0}^s}\sum_{G\in A_{\clG}} |G'|^{s-\veps} |G'|^{\veps} < \frac{1}{\kappa^{\veps j_0}B_{j_0}^s} \sum_{G\in A_{\clG}} |G'|^{s-\veps} \nonumber\\
&< \sum_{G\in A_{\clG}} |G'|^{s-\veps}.\nonumber
\end{align*}
Take $\clG':=\{G':G\in\clG\}$, then \eqref{Ec1GJL04} applied to $\mathfrak{A}=\{A\in\clA: A_{\clG}\neq\vac\}$ yields
\[
c<\Lambda_s(\mathfrak{A})\leq \Lambda_{s-\veps}(\clG')\leq \Lambda_{s-\veps}(\clG).
\]
Since $\veps>0$ and $\clG$ were arbitrary, $\dim_H \bfA_{\infty}\geq s$.
\end{proof}

In practice, rather than \eqref{Ec1GJL04}, we will use a stronger condition. Namely, keeping the notation of Lemma \ref{LeGJL01},  
\begin{equation}\label{EcW1GJL}
\forall n\in\Na \quad \forall A\in\clA_n\quad
\sum_{B\in D(A)} |B|^s\geq |A|^s.
\end{equation}

Let us check that \eqref{EcW1GJL} indeed implies \eqref{Ec1GJL04}. Take a finite cover $\mfX\subseteq \clA$ of $\mathbf{A}_{\infty}$. We may assume that the elements of $\mfX$ are disjoint by pairs. Let $n$ be the maximal integer satisfying $\clA_n\cap \mfX\neq \vac$. Take $B\in \clA_n$ and $A\in\clA_{n-1}$ such that $B\in D(A)$. By pairwise disjointness and since $\mfX$ covers $\bfA_{\infty}$, $D(A)\subseteq \mfX$. Thus, we may replace the members of $D(A)$ by $A$ to obtain a new covering $\mfX'$ satisfying $\Lambda_s (\mfX)\geq \Lambda_s (\mfX')$. Repeating the process, we eventually arrive at $\clA_0$ and \eqref{Ec1GJL04} follows.

\begin{lem01}[Second Generalized Jarn\'ik Lemma]\label{LeGJL02}
Let $\clA$ be a  family of compact sets satisfying conditions i., iii., iv., v. of Definition \ref{DefDSTL} and such that the each $\clA_k$, $k\geq 1$, is at most countable.  

Let $s>0$. If 
\begin{equation}\label{EcLeGJL02}
\forall k\in\Na \quad \forall A\in\clA_k \quad
\sum_{B\in D(A)} |B|^s  \leq |A|^s,
\end{equation}
then $\dim_H A_{\infty}\leq s$.
\end{lem01}
\begin{proof}
We keep the statement's notation. Note that $(\Lambda_s(\clA_j))_{j\geq 1}$ is decreasing, because
\[
\forall j\in \Na_0 \quad 
\Lambda_s(\clA_{j+1}) = \sum_{B\in\clA_{j+1}} |B|^s =\sum_{A\in\clA_{j}} \sum_{B\in D(A)} |B|^s \leq \sum_{A\in \clA_{j}} |A|^s\leq \Lambda_s(\clA_{j}).
\]
Since every $\clA_j$ covers $\bfA_{\infty}$ and $d_j(\clA)\to 0$ as $j\to\infty$, $\dim_H \bfA_{\infty}\leq s$.
\end{proof}
\begin{rema1000}
The conclusion of Lemma \ref{LeGJL01} (resp. Lemma \ref{LeGJL02}) is true if \eqref{EcW1GJL} (resp. \eqref{EcLeGJL02}) only holds for every sufficiently large $n\in\Na$.
\end{rema1000}
\section{Preliminary Lemmas}

Define the following sets
\begin{align*}
E&:=\{ z=[0;a_1,a_2,a_3,\ldots]\in\mfF: \lim_{n\to\infty} |a_n|=\infty\} \nonumber\\
\forall L,M>0 \quad E_L^M &:= \left\{ z=[0;a_1,a_2,a_3,\ldots] \in\mfF : \forall n\in\Na \quad L \leq |a_n|\leq M \right\}.
\end{align*}

\begin{lem01}\label{Le4.1}
\begin{enumerate}[1.]
\item For every $L\geq \sqrt{8}$ we have $1\leq \displaystyle\lim_{M\to\infty}\,\dim_H E_L^M$.
\item $\displaystyle\lim_{L\to\infty} \dim_H E_L\leq 1$.
\end{enumerate}
\end{lem01}

For any $z\in\Cx$ write $\|z\|=\max\{|\Re z|,|\Im z|\}$. If $f,g:\Za[i]\to\RE_{>0}$ satisfy $f(n)\leq g(n)$ for all $n\in\Na$, we define 
\[
E_{f,g} :=\left\{ [0;a_1,a_2,a_3,\ldots] \in\mfF : \forall n\in\Na \quad f(n)\leq \|a_n\|\leq  g(n) \right\}.
\]
Note that for every $z\in\Cx$ we have $\|z\| \leq |z|\leq \sqrt{2}\|z\|$.
\begin{lem01}\label{Le4.3}
Let $f,g:\Na\to\RE_{>0}$ be functions such that for some fixed $c'$ with $0<c'<1$ we have $\sqrt{8}\leq f(n)<c'g(n)$ for all $n\in\Na$, $f(n)\to\infty$ when $n\to\infty$, and 
\[
\limsup_{n\to\infty} \frac{\log \log g(n) }{\log n}<1. 
\]
Then, $\dim_{H} E_{f,g} = 1$.
\end{lem01}

\section{Proofs of Main Results}
\begin{proof}[Proof of Theorem \ref{TeoGero}]
It is enough to prove the theorem on $\mfF$, because the zeroth term of a HCF expansion does not impose any restrictions on the subsequent terms. By $E \subseteq E_{L}'$ for $L>\sqrt{8}$, the second part of Lemma \ref{Le4.1}, and Lemma \ref{LeDesc01},
\[
\dim_H E \leq 1. 
\]
Let $f$ and $g$ be as in Lemma \ref{Le4.3}, then $E_{f,g} \subseteq E$ and
\[
1
\leq 
\dim_H E_{f,g}
\leq 
\dim_H E. 
\]
\end{proof}
\begin{rema1000}
If $f:\Na\to\RE_{>0}$, $f(n)\to\infty$ as $n\to\infty$ and 
\[
\limsup_{n\to\infty} \frac{\log(\log f(n))}{\log n} <1, 
\]
then
\[
E_f=\left\{ z=[0;a_1,a_2,\ldots]: \forall n\in\Na \quad f(n)\leq \|a_n\|\right\}
\]
has Hausdorff dimension $1$, because we can easily find an adequate $g$ such that $E_{f,g}\subseteq E_f \subseteq E$ and $\dim_H E_{f,g} = 1$.
\end{rema1000}

\begin{proof}[Proof of Corollary \ref{CoroGero}]
Let $\veps>0$. By Lemma \ref{Le4.1}, there are $L=L(\veps)\in\Na$ and $M=M(L,\veps)\in\Na$ such that
\[
1-\veps \leq \dim_{H} E_L^M \leq \dim_{H} E_L \leq 1+ \veps.
\]
Since an irrational complex number belongs to $\bad_{\Cx}$ if and only if its HCF is bounded, $E_L^M\subseteq E_L\cap \bad_{\Cx}\subseteq E_L$ and the result follows.
\end{proof}

\section{Proofs of Preliminary Lemmas}
\subsection{Proof of Lemma \ref{Le4.1}}
\begin{proof}[Proof of Lemma \ref{Le4.1}]
\begin{enumerate}[i.]
\item The result will follow from Lemma \ref{LeGJL01} applied to the family of compact sets $\clA=\bigcup_j \clA_j$, where
\begin{equation}\label{EcL4101}
\forall j\in\Na \quad 
\clA_j :=\left\{ \Cl\left( \clC_j(\bfa ) \cap E_{\sqrt{8}} \right) :  \bfa\in\Omega_0^{\HCF}, \quad \forall n\in\{1,\ldots,j\} \quad L \leq |a_n|\leq M \right\}
\end{equation}
for any given $L\geq\sqrt{8}$ and large enough $M$.

Let $L\geq \sqrt{8}$ and let $M>L$ be such that \eqref{EcL4101} has more than $1$ element. Take $\bfa\in E_{\sqrt{8}}$. By \eqref{Ec2.2} and Proposition \ref{Propo2.1} \ref{Propo2.1_iii}., there is a constant $\psi>1$ such that
\[
\forall n\in\Na \quad \left| \clC_n(\bfa)\cap E_{\sqrt{8}}\right|\leq \frac{2}{\psi^{n-1}}.
\]
Since $\bfa$ is arbitrary, \eqref{Ec1GJL01} follows.

After taking $k>0$ as in Lemma \ref{Le2.0}, we consider the sequence $B_n=\tfrac{k}{2}$ for all $n\in\Na$. Clearly, \eqref{Ec1GJL02} holds for $(B_n)_{n\geq 1}$ and, by Proposition \ref{Propo2.1} part \ref{Propo2.1_iv} and Lemmas \ref{Le2.0} and \ref{Le2.1}, $(B_n)_{n\geq 1}$ also satisfies \eqref{Ec1GJL03}. (The proof of Lemma \ref{Le4.3} contains a similar argument done in full detail.)

In order to verify \eqref{EcW1GJL}, we further assume that $M$ is large enough for 
\[
\#\{b\in\Za[i]: L\leq |b|\leq M\} \geq M^2.
\]
to hold. Take $n\in\Na$ and let $\bfb=(b_1,\ldots,b_n)\in\Za[i]^n$ be such that $L\leq |b_j|\leq M$ for every $j$ and denote its $\clQ$-pair by $\sepn,\seqn$. Consider $\gamma$ as in Lemma \ref{Le2.1}; then, for every $s>0$
\begin{align*}
\sum_{L\leq |c|\leq M} \left| \clC_{n+1}(\bfb c)\cap E_{\sqrt{8}}\right|^{s} &\geq 
\frac{\gamma^{s}}{(M+1)^{2s}} \sum_{L\leq |c|\leq M } \frac{1}{|q_n|^{2s}} \geq \frac{\gamma^{s} M^2 }{(M+1)^{2s}} \frac{1}{|q_n|^{2s}} \nonumber\\
&\geq \frac{ M^{2-2s} \gamma^{s}} { \left(1+ \frac{2}{M} \right)^{2s} 2^{s}}|\clC_n(\bfb)\cap E_{\sqrt{8}}|^{s}. \nonumber
\end{align*}
The coefficient of $|\clC_n(\bfb)\cap E_{\sqrt{8}}|^{s}$ is at least $1$ if and only if
\[
(2-2s) \log M + s\log\frac{\gamma}{2} - 2s\log\left( 1+\frac{1}{M}\right) \geq 0,
\]
which, after rearranging the terms, is equivalent to
\[
s \leq \frac{2 \log M }{2\log M + 2\log\left( 1+\frac{1}{M}\right) - \log\frac{\gamma}{2}}<1.
\]
Since the middle term tends to $1$ when $M\to\infty$, \eqref{EcW1GJL} holds for any given $0<s<1$ as long as $M=M(s,L)$ is large enough. For such $s,M$, Lemma \ref{LeGJL01} implies $s\leq \dim_H E_L^M$.  
\item The result follows from Lemma \ref{LeGJL02} applied to the family of compact sets $\clA^L= \bigcup_j \clA^L_j$, where
\[
\forall n\in\Na \quad
\clA^L_n:=\left\{ \Cl\left(\clC_n(\bfa)\cap E_{\sqrt{8}} \right): \bfa=(a_1,\ldots,a_n)\in\Za[i]^n, \; L\leq \min_j |a_j| \right\},
\]
and sufficiently large $L>2$.

We shall only verify \eqref{EcLeGJL02}, the other conditions are obtained as above. Consider $L>0$ and $0<\veps<1$. Take $n\in\Na$, $\bfa=(a_1,\ldots,a_n)\in\Za[i]^n$ with $|a_j|\geq \sqrt{8}$ for every $j$ and with $\clQ$-pair $\sepn, \seqn$. Then, by Lemma  \ref{Le2.1} and the recurrence defining $\seqn$,
\[
\sum_{|b|\geq L} |\clC_{n+1}(\bfa\, b)\cap E_{\sqrt{8}}|^{1+\veps} 
\leq \frac{2^{1+\veps}}{|q_n|^{2+2\veps}} \sum_{|b|\geq L } \frac{1}{\left|b + \frac{q_{n-1}}{q_n}\right|^{2+2\veps}} 
\leq
\frac{2^{3+3\veps}}{|q_n|^{2+2\veps}} \sum_{|b|\geq L } \frac{1}{|b|^{2+2\veps}}.
\]
In the last inequality we used $L>2$. Hence, for some constants $c_1=c_1(\veps)>0$, $c_2=c_2(\veps)>0$
\[
\sum_{|b|\geq L} |\clC_{n+1}(\bfa\, b)\cap E_{\sqrt{8}}|^{1+\veps} 
\leq
\frac{c_1}{|q_n|^{2+2\veps}} \int_L^{\infty}\frac{\md r}{r^{1+2\veps}}
\leq
\frac{c_2}{L^{2\veps}} |\clC_n(\bfa)\cap E_{\sqrt{8}}|^{1+\veps}.
\]
By the Second Generalized Jarn\'ik Lemma, for large $L=L(\veps)$ we have 
\[
\dim_H E_L\leq 1+\veps. 
\]
\end{enumerate}
\end{proof}

\subsection{Proof of Lemma \ref{Le4.3}}
\begin{proof}[Proof of Lemma \ref{Le4.3}]
Keep the statement's notation. By $E_{f,g}\subseteq E$ and Lemma \ref{Le4.1}, $\dim_H E_{f,g}\leq 1$. The inequality $\dim_H E_{f,g}\geq 1$ will follow from Lemma \ref{LeGJL01} applied to the family of compact sets
\[
\left\{ \Cl\left( \clC_n(\bfa)\cap E_{\sqrt{8}}\right): n\in\Na, \; \bfa=(a_1,\ldots,a_n)\in\Za[i]^n, \; \forall j\in\{1,\ldots,n\}\quad f(j)\leq \|a_j\| \leq g(j)\right\}.
\]
The condition \eqref{Ec1GJL01} clearly holds by part \ref{Propo2.1_iii}. of Proposition \ref{Propo2.1} and Lemma \ref{Le2.1}.

Take $n\in\Na$. Let $\bfa\in\Za[i]^n$ be such that
\[
\forall j\in\{1,\ldots,n\} \quad
f(j)\leq \|a_j\|\leq g(j),
\]
and take $b,c\in\Za[i]$ satisfying $f(n+1)\leq \|b \|, \|c \| \leq g(n+1)$, and $b\neq c$. Note that the $\clQ$-pairs of $\bfa b$ and $\bfa c$ coincide except for the last term. Denote by $(p_j)_{j=0}^{n+1}$, $(q_j)_{j=0}^{n+1}$ the $\clQ$-pair of $\bfa b$ and let $p_{n+1}'$, $q_{n+1}'$ be the last term of the $\clQ$-pair of $\bfa c$. Take $k$ as in Lemma \ref{Le2.0}. Then, there exists an absolute constant $c_3>0$ such that for every $\alpha,\beta\in E_{\sqrt{8}}$
\begin{align*}
d\left(\clC_{n+1}(\bfa\,b)\cap E_{\sqrt{8}}, \clC_{n+1}(\bfa\,c)\cap E_{\sqrt{8}}\right)&\geq
\left| \frac{\beta p_n + p_{n+1}}{\beta q_n + q_{n+1}} - \frac{\alpha p_n + p_{n+1}'}{\alpha q_n + q_{n+1}'} \right| \nonumber\\
&= 
\frac{1}{|q_n|^2}  \frac{\left| c+\alpha - (b+\beta)\right|}{\left|\left(b+z + \frac{q_{n-1}}{q_n}\right)\left(c+w+\frac{q_{n-1}}{q_n}\right)\right|}  \nonumber\\
&\geq 
\frac{1}{|q_n|^2} \frac{k}{(g(n)+2)^2} 
\geq 
\frac{c_3}{g(n)^2} \left|\clC_{n+1}(\bfa)\cap E_{\sqrt{8}}\right|.\nonumber
\end{align*}
Since $n$ and $\bfa$ are arbitrary, the sequence $(B_n)_{n\geq 1}$ given by $B_n = c_3g(n)^{-2}$,  $n\in\Na$, satisfies \eqref{Ec1GJL02} and \eqref{Ec1GJL03}.

It remains to show that \eqref{EcW1GJL} holds. Take $s>0$. Using repeatedly Lemma \ref{Le2.1} and taking $\gamma$ as in \eqref{Ec2.2}, there are absolute constants $c_4,c_5>0$ for which
\begin{align*}
\sum_{f(n+1)\leq \|b\|\leq g(n+1)} |\clC_{n+1}(\bfa b)\cap E_{\sqrt{8}}|^{s} &\geq \frac{\gamma^{s}}{|q_n|^{2s}} \sum_{f(n+1)\leq \|b\|\leq g(n+1)} \left| b + \frac{q_{n-1}}{q_{n}} \right|^{-2s} \nonumber\\
&\geq c_4\frac{\gamma^{s}}{|q_n|^{2s}}  \sum_{f(n+1)\leq \|b\|\leq g(n+1)} (g(n+1)+1)^{-2s}\nonumber\\
&\geq  c_4\frac{\gamma^{s}}{|q_n|^{2s}} \frac{g(n+1)^2-f(n+1)^2}{(g(n+1)+1)^{2s}} \nonumber\\
&\geq c_5 \frac{\gamma^{s}}{2^{s}} \frac{\left(1-\frac{f(n+1)^2}{g(n+1)^2}\right)}{\left(1+ \frac{1}{g(n+1)}\right)^{2s}} g(n+1)^{2-2s}|\clC_n(\bfa)\cap E_{\sqrt{8}}|^{s}. \nonumber
\end{align*}
After taking logarithms and rearranging terms, we see that the coefficient of $|\clC_n(\bfa)\cap E_{\sqrt{8}}|^{s}$ is at least $1$ if and only if
\begin{equation}\label{EcUlt}
\frac{2\log(g(n+1)) + \log\left( 1- \frac{f(n+1)^2}{g(n+1)^2}\right) + \log(c_5)}{2\log(g(n+1)) + 2\log\left( 1 + \frac{1}{g(n+1)}\right) - \log\left(\frac{\gamma}{2}\right)}>s.
\end{equation}

Since for some $0<c'<1$ and every $j$ the inequality $f(j)\leq c' g(j)$ is true, the left hand side of \eqref{EcUlt} tends to $1$ as $n$ tends to $\infty$. Therefore, \eqref{EcUlt} holds for any $0<s<1$ and all $n\in\Na$ larger than some $N=N(s)$. Finally, by Lemma \ref{LeGJL01}, $\dim_H E_{f,g}\geq 1$.
\end{proof}

\begin{flushleft}
\textbf{Acknowledgements}
\end{flushleft}

This research was supported by CONACyT, Mexico, grant 410695.

\begin{flushleft}
Gerardo Gonz\'alez Robert \par
Facultad de Ciencias, \par
Universidad Nacional Aut\'onoma de M\'exico, \par
Circuito Exterior S/N, C.U., Coyoac\'an, 04510 \par
Mexico City, Mexico \par
\texttt{gerardogonrob@ciencias.unam.mx}
\end{flushleft}
\end{document}